\documentclass{amsart}
\usepackage{amsfonts}

\newtheorem{theorem}{Theorem}
\theoremstyle{plain}

\newtheorem{corollary}{Corollary}

\newtheorem{lemma}{Lemma}

\newtheorem{remark}{Remark}

\numberwithin{equation}{section}
\input{tcilatex}

\begin{document}
\title[Pe\v{c}ari\'{c} Inequality]{On Pe\v{c}ari\'{c}'s Inequality in Inner Product Spaces}
\author{S.S. Dragomir}
\address{School of Computer Science and Mathematics\\
Victoria University of Technology\\
PO Box 14428, MCMC \\
Victoria 8001, Australia.}
\email{sever.dragomir@vu.edu.au}
\urladdr{http://rgmia.vu.edu.au/SSDragomirWeb.html}
\date{12 June, 2003.}

\begin{abstract}
Some related results to Pe\v{c}ari\'{c}'s inequality in inner product spaces
that generalises Bombieri's inequality, are given.
\end{abstract}

\maketitle

\section{Introduction}

In 1992, J.E. Pe\v{c}ari\'{c} \cite{3b} proved the following inequality for
vectors in complex inner product spaces $\left( H;\left( \cdot ,\cdot
\right) \right) $.

\begin{theorem}
\label{t1.1}Suppose that $x,y_{1},\dots ,y_{n}$ are vectors in $H$ and $%
c_{1},\dots ,c_{n}$ are complex numbers. Then the following inequalities 
\begin{eqnarray}
\left| \sum\limits_{i=1}^{n}c_{i}\left( x,y_{i}\right) \right| ^{2} &\leq
&\left\| x\right\| ^{2}\sum\limits_{i=1}^{n}\left| c_{i}\right| ^{2}\left(
\sum\limits_{j=1}^{n}\left| \left( y_{i},y_{j}\right) \right| \right)
\label{1.1} \\
&\leq &\left\| x\right\| ^{2}\sum\limits_{i=1}^{n}\left| c_{i}\right|
^{2}\max_{1\leq i\leq n}\left( \sum\limits_{j=1}^{n}\left| \left(
y_{i},y_{j}\right) \right| \right) ,  \notag
\end{eqnarray}
hold.
\end{theorem}

He also showed that for $c_{i}=\overline{\left( x,y_{i}\right) },i\in
\left\{ 1,...,n\right\} ,$ one gets 
\begin{align}
\left( \sum\limits_{i=1}^{n}\left| \left( x,y_{i}\right) \right| ^{2}\right)
^{2}& \leq \left\| x\right\| ^{2}\sum\limits_{i=1}^{n}\left| \left(
x,y_{i}\right) \right| ^{2}\left( \sum\limits_{j=1}^{n}\left| \left(
y_{i},y_{j}\right) \right| \right)  \label{1.2} \\
& \leq \left\| x\right\| ^{2}\sum\limits_{i=1}^{n}\left| \left(
x,y_{i}\right) \right| ^{2}\max_{1\leq i\leq n}\left(
\sum\limits_{j=1}^{n}\left| \left( y_{i},y_{j}\right) \right| \right) , 
\notag
\end{align}
which improves Bombieri's result \cite{1b} (see also \cite[p. 394]{2b}) 
\begin{equation}
\sum\limits_{i=1}^{n}\left| \left( x,y_{i}\right) \right| ^{2}\leq \left\|
x\right\| ^{2}\max_{1\leq i\leq n}\left( \sum\limits_{j=1}^{n}\left| \left(
y_{i},y_{j}\right) \right| \right) .  \label{1.3}
\end{equation}
Note that (\ref{1.3}) is in its turn a natural generalisation of \textit{%
Bessel's inequality} 
\begin{equation}
\sum\limits_{i=1}^{n}\left| \left( x,e_{i}\right) \right| ^{2}\leq \left\|
x\right\| ^{2},\ \ x\in H,  \label{1.4}
\end{equation}
which holds for the orthornormal vectors $\left( e_{i}\right) _{1\leq i\leq
n}.$

In this paper we point out some related results to Pe\v{c}ari\'{c}'s
inequality (\ref{1.1}). Some results of Bombieri type are also mentioned.

\section{Preliminary Results}

We start with the following lemma that is interesting in its own right.

\begin{lemma}
\label{l2.1}Let $z_{1},\dots ,z_{n}\in H$ and $\alpha _{1},\dots ,\alpha
_{n}\in \mathbb{K}$. Then one has the inequalities: 
\begin{equation}
\left\| \sum_{i=1}^{n}\alpha _{i}z_{i}\right\| ^{2}\leq \left(
\sum\limits_{i=1}^{n}\left| \alpha _{i}\right| ^{p}\left(
\sum\limits_{j=1}^{n}\left| \left( z_{i},z_{j}\right) \right| \right)
\right) ^{\frac{1}{p}}\left( \sum\limits_{i=1}^{n}\left| \alpha _{i}\right|
^{q}\left( \sum\limits_{j=1}^{n}\left| \left( z_{i},z_{j}\right) \right|
\right) \right) ^{\frac{1}{q}}  \label{2.1}
\end{equation}
\begin{equation*}
\leq \left\{ 
\begin{array}{l}
\max\limits_{1\leq i\leq n}\left| \alpha _{i}\right|
^{2}\sum\limits_{i,j=1}^{n}\left| \left( z_{i},z_{j}\right) \right| ; \\ 
\\ 
\max\limits_{1\leq i\leq n}\left| \alpha _{i}\right| \left(
\sum\limits_{i=1}^{n}\left| \alpha _{i}\right| ^{\gamma q}\right) ^{\frac{1}{%
\gamma q}}\left( \sum\limits_{i,j=1}^{n}\left| \left( z_{i},z_{j}\right)
\right| \right) ^{\frac{1}{p}}\left( \sum\limits_{i=1}^{n}\left(
\sum\limits_{j=1}^{n}\left| \left( z_{i},z_{j}\right) \right| \right)
^{\delta }\right) ^{\frac{1}{\delta q}}, \\ 
\hfill \ \ \text{if}\ \gamma >1,\ \frac{1}{\gamma }+\frac{1}{\delta }=1; \\ 
\\ 
\max\limits_{1\leq i\leq n}\left| \alpha _{i}\right| \left(
\sum\limits_{i=1}^{n}\left| \alpha _{i}\right| ^{q}\right) ^{\frac{1}{q}%
}\left( \sum\limits_{i,j=1}^{n}\left| \left( z_{i},z_{j}\right) \right|
\right) ^{\frac{1}{p}}\max\limits_{1\leq i\leq n}\left(
\sum\limits_{j=1}^{n}\left| \left( z_{i},z_{j}\right) \right| \right) ^{%
\frac{1}{q}}; \\ 
\\ 
\max\limits_{1\leq i\leq n}\left| \alpha _{i}\right| \left(
\sum\limits_{i=1}^{n}\left| \alpha _{i}\right| ^{\alpha p}\right) ^{\frac{1}{%
\alpha p}}\left( \sum\limits_{i,j=1}^{n}\left| \left( z_{i},z_{j}\right)
\right| \right) ^{\frac{1}{q}}\left( \sum\limits_{i=1}^{n}\left(
\sum\limits_{j=1}^{n}\left| \left( z_{i},z_{j}\right) \right| \right)
^{\beta }\right) ^{\frac{1}{\beta q}}, \\ 
\hfill \ \ \text{if}\ \alpha >1,\ \frac{1}{\alpha }+\frac{1}{\beta }=1; \\ 
\left( \sum\limits_{i=1}^{n}\left| \alpha _{i}\right| ^{\alpha p}\right) ^{%
\frac{1}{\alpha p}}\left( \sum\limits_{i=1}^{n}\left| \alpha _{i}\right|
^{\gamma q}\right) ^{\frac{1}{\gamma q}}\left( \sum\limits_{i=1}^{n}\left(
\sum\limits_{j=1}^{n}\left| \left( z_{i},z_{j}\right) \right| \right)
^{\beta }\right) ^{\frac{1}{p\beta }} \\ 
\ \ \ \ \ \times \left( \sum\limits_{i=1}^{n}\left(
\sum\limits_{j=1}^{n}\left| \left( z_{i},z_{j}\right) \right| \right)
^{\delta }\right) ^{\frac{1}{\delta q}}\hfill \ \text{if}\ \alpha >1,\ \frac{%
1}{\alpha }+\frac{1}{\beta }=1\text{ and }\ \gamma >1,\ \frac{1}{\gamma }+%
\frac{1}{\delta }=1; \\ 
\\ 
\left( \sum\limits_{i=1}^{n}\left| \alpha _{i}\right| ^{q}\right) ^{\frac{1}{%
q}}\left( \sum\limits_{i=1}^{n}\left| \alpha _{i}\right| ^{\alpha p}\right)
^{\frac{1}{\alpha p}}\max\limits_{1\leq i\leq n}\left(
\sum\limits_{j=1}^{n}\left| \left( z_{i},z_{j}\right) \right| \right) ^{%
\frac{1}{q}}\left( \sum\limits_{i=1}^{n}\left( \sum\limits_{j=1}^{n}\left|
\left( z_{i},z_{j}\right) \right| \right) ^{\beta }\right) ^{\frac{1}{p\beta 
}}, \\ 
\hfill \ \ \text{if}\ \alpha >1,\ \frac{1}{\alpha }+\frac{1}{\beta }=1; \\ 
\\ 
\max\limits_{1\leq i\leq n}\left| \alpha _{i}\right| \left(
\sum\limits_{i=1}^{n}\left| \alpha _{i}\right| ^{p}\right) ^{\frac{1}{p}%
}\max\limits_{1\leq i\leq n}\left( \sum\limits_{j=1}^{n}\left| \left(
z_{i},z_{j}\right) \right| \right) ^{\frac{1}{p}}\left(
\sum\limits_{i,j=1}^{n}\left| \left( z_{i},z_{j}\right) \right| \right) ^{%
\frac{1}{q}}; \\ 
\\ 
\left( \sum\limits_{i=1}^{n}\left| \alpha _{i}\right| ^{p}\right) ^{\frac{1}{%
p}}\left( \sum\limits_{i=1}^{n}\left| \alpha _{i}\right| ^{\gamma q}\right)
^{\frac{1}{\gamma q}}\max\limits_{1\leq i\leq n}\left(
\sum\limits_{j=1}^{n}\left| \left( z_{i},z_{j}\right) \right| \right) ^{%
\frac{1}{p}}\left( \sum\limits_{i=1}^{n}\left( \sum\limits_{j=1}^{n}\left|
\left( z_{i},z_{j}\right) \right| \right) ^{\delta }\right) ^{\frac{1}{%
\delta q}}, \\ 
\hfill \ \ \text{if}\ \gamma >1,\ \frac{1}{\gamma }+\frac{1}{\delta }=1; \\ 
\left( \sum\limits_{i=1}^{n}\left| \alpha _{i}\right| ^{p}\right) ^{\frac{1}{%
p}}\left( \sum\limits_{i=1}^{n}\left| \alpha _{i}\right| ^{q}\right) ^{\frac{%
1}{q}}\max\limits_{1\leq i\leq n}\left( \sum\limits_{j=1}^{n}\left| \left(
z_{i},z_{j}\right) \right| \right) ,
\end{array}
\right.
\end{equation*}
where $p>1,$ $\frac{1}{p}+\frac{1}{q}=1.$
\end{lemma}

\begin{proof}
We observe that 
\begin{eqnarray}
\left\| \sum_{i=1}^{n}\alpha _{i}z_{i}\right\| ^{2} &=&\left(
\sum_{i=1}^{n}\alpha _{i}z_{i},\sum_{j=1}^{n}\alpha _{j}z_{j}\right)
\label{2.2} \\
&=&\sum_{i=1}^{n}\sum_{j=1}^{n}\alpha _{i}\overline{\alpha _{j}}\left(
z_{i},z_{j}\right) =\left| \sum_{i=1}^{n}\sum_{j=1}^{n}\alpha _{i}\overline{%
\alpha _{j}}\left( z_{i},z_{j}\right) \right|  \notag \\
&\leq &\sum_{i=1}^{n}\sum_{j=1}^{n}\left| \alpha _{i}\right| \left| \alpha
_{j}\right| \left| \left( z_{i},z_{j}\right) \right| =:M.  \notag
\end{eqnarray}
If one uses the H\"{o}lder inequality for double sums, i.e., we recall it 
\begin{equation}
\sum\limits_{i,j=1}^{n}m_{ij}a_{ij}b_{ij}\leq \left(
\sum\limits_{i,j=1}^{n}m_{ij}a_{ij}^{p}\right) ^{\frac{1}{p}}\left(
\sum\limits_{i,j=1}^{n}m_{ij}b_{ij}^{q}\right) ^{\frac{1}{q}},  \label{2.3}
\end{equation}
where $m_{ij},a_{ij},b_{ij}\geq 0,$ $\frac{1}{p}+\frac{1}{q}=1,$ $p>1;$ then 
\begin{align}
M& \leq \left( \sum\limits_{i,j=1}^{n}\left| \left( z_{i},z_{j}\right)
\right| \left| \alpha _{i}\right| ^{p}\right) ^{\frac{1}{p}}\left(
\sum\limits_{i,j=1}^{n}\left| \left( z_{i},z_{j}\right) \right| \left|
\alpha _{i}\right| ^{q}\right) ^{\frac{1}{q}}  \label{2.4} \\
& =\left( \sum_{i=1}^{n}\left| \alpha _{i}\right| ^{p}\left(
\sum_{j=1}^{n}\left| \left( z_{i},z_{j}\right) \right| \right) \right) ^{%
\frac{1}{p}}\left( \sum_{i=1}^{n}\left| \alpha _{i}\right| ^{q}\left(
\sum_{j=1}^{n}\left| \left( z_{i},z_{j}\right) \right| \right) \right) ^{%
\frac{1}{q}},  \notag
\end{align}
and the first inequality in (\ref{2.1}) is proved.

Observe that 
\begin{equation*}
\sum_{i=1}^{n}\left| \alpha _{i}\right| ^{p}\left( \sum_{j=1}^{n}\left|
\left( z_{i},z_{j}\right) \right| \right) \leq \left\{ 
\begin{array}{l}
\max\limits_{1\leq i\leq n}\left| \alpha _{i}\right|
^{p}\sum\limits_{i,j=1}^{n}\left| \left( z_{i},z_{j}\right) \right| ; \\ 
\\ 
\left( \sum\limits_{i=1}^{n}\left| \alpha _{i}\right| ^{\alpha p}\right) ^{%
\frac{1}{\alpha }}\left( \sum\limits_{i=1}^{n}\left(
\sum\limits_{j=1}^{n}\left| \left( z_{i},z_{j}\right) \right| \right)
^{\beta }\right) ^{\frac{1}{\beta }} \\ 
\hfill \ \text{if}\ \alpha >1,\ \frac{1}{\alpha }+\frac{1}{\beta }=1; \\ 
\\ 
\sum\limits_{i=1}^{n}\left| \alpha _{i}\right| ^{p}\max\limits_{1\leq i\leq
n}\left( \sum\limits_{j=1}^{n}\left| \left( z_{i},z_{j}\right) \right|
\right) ;
\end{array}
\right.
\end{equation*}
giving 
\begin{multline}
\left( \sum_{i=1}^{n}\left| \alpha _{i}\right| ^{p}\left(
\sum_{j=1}^{n}\left| \left( z_{i},z_{j}\right) \right| \right) \right) ^{%
\frac{1}{p}}  \label{2.5} \\
\leq \left\{ 
\begin{array}{l}
\max\limits_{1\leq i\leq n}\left| \alpha _{i}\right| \left(
\sum\limits_{i,j=1}^{n}\left| \left( z_{i},z_{j}\right) \right| \right) ^{%
\frac{1}{p}};\hfill \\ 
\\ 
\left( \sum\limits_{i=1}^{n}\left| \alpha _{i}\right| ^{\alpha p}\right) ^{%
\frac{1}{\alpha p}}\left( \sum\limits_{i=1}^{n}\left(
\sum\limits_{j=1}^{n}\left| \left( z_{i},z_{j}\right) \right| \right)
^{\beta }\right) ^{\frac{1}{\beta p}}\hfill \ \text{if}\ \alpha >1,\ \frac{1%
}{\alpha }+\frac{1}{\beta }=1; \\ 
\\ 
\left( \sum\limits_{i=1}^{n}\left| \alpha _{i}\right| ^{p}\right) ^{\frac{1}{%
p}}\max\limits_{1\leq i\leq n}\left( \sum\limits_{j=1}^{n}\left| \left(
z_{i},z_{j}\right) \right| \right) ^{\frac{1}{p}}.
\end{array}
\right.
\end{multline}
Similarly, we have 
\begin{multline}
\left( \sum_{i=1}^{n}\left| \alpha _{i}\right| ^{q}\left(
\sum_{j=1}^{n}\left| \left( z_{i},z_{j}\right) \right| \right) \right) ^{%
\frac{1}{q}}  \label{2.6} \\
\leq \left\{ 
\begin{array}{l}
\max\limits_{1\leq i\leq n}\left| \alpha _{i}\right| \left(
\sum\limits_{i,j=1}^{n}\left| \left( z_{i},z_{j}\right) \right| \right) ^{%
\frac{1}{q}}\hfill \\ 
\\ 
\left( \sum\limits_{i=1}^{n}\left| \alpha _{i}\right| ^{\gamma q}\right) ^{%
\frac{1}{\gamma q}}\left( \sum\limits_{i=1}^{n}\left(
\sum\limits_{j=1}^{n}\left| \left( z_{i},z_{j}\right) \right| \right)
^{\delta }\right) ^{\frac{1}{\delta q}}\hfill \ \text{if}\ \gamma >1,\ \frac{%
1}{\gamma }+\frac{1}{\delta }=1; \\ 
\\ 
\left( \sum\limits_{i=1}^{n}\left| \alpha _{i}\right| ^{q}\right) ^{\frac{1}{%
q}}\max\limits_{1\leq i\leq n}\left( \sum\limits_{j=1}^{n}\left| \left(
z_{i},z_{j}\right) \right| \right) ^{\frac{1}{q}}.
\end{array}
\right.
\end{multline}
Using (\ref{2.1}) and (\ref{2.5}) -- (\ref{2.6}), we deduce the 9
inequalities in the second part of (\ref{2.2}).
\end{proof}

If we choose $p=q=2,$ then the following result holds.

\begin{corollary}
\label{c2.2}If $z_{1},\dots ,z_{n}\in H$ and $\alpha _{1},\dots ,\alpha
_{n}\in \mathbb{K}$, then one has 
\begin{equation}
\left\| \sum_{i=1}^{n}\alpha _{i}z_{i}\right\| ^{2}\leq
\sum\limits_{i=1}^{n}\left| \alpha _{i}\right| ^{2}\left(
\sum\limits_{j=1}^{n}\left| \left( z_{i},z_{j}\right) \right| \right)
\label{2.7}
\end{equation}
\begin{equation*}
\leq \left\{ 
\begin{array}{l}
\max\limits_{1\leq i\leq n}\left| \alpha _{i}\right|
^{2}\sum\limits_{i,j=1}^{n}\left| \left( z_{i},z_{j}\right) \right| ; \\ 
\\ 
\max\limits_{1\leq i\leq n}\left| \alpha _{i}\right| \left(
\sum\limits_{i=1}^{n}\left| \alpha _{i}\right| ^{2\gamma }\right) ^{\frac{1}{%
2\gamma }}\left( \sum\limits_{i,j=1}^{n}\left| \left( z_{i},z_{j}\right)
\right| \right) ^{\frac{1}{2}}\left( \sum\limits_{i=1}^{n}\left(
\sum\limits_{j=1}^{n}\left| \left( z_{i},z_{j}\right) \right| \right)
^{\delta }\right) ^{\frac{1}{2\delta }}, \\ 
\hfill \ \ \text{if}\ \gamma >1,\ \frac{1}{\gamma }+\frac{1}{\delta }=1; \\ 
\max\limits_{1\leq i\leq n}\left| \alpha _{i}\right| \left(
\sum\limits_{i=1}^{n}\left| \alpha _{i}\right| ^{2}\right) ^{\frac{1}{2}%
}\left( \sum\limits_{i,j=1}^{n}\left| \left( z_{i},z_{j}\right) \right|
\right) ^{\frac{1}{2}}\max\limits_{1\leq i\leq n}\left(
\sum\limits_{j=1}^{n}\left| \left( z_{i},z_{j}\right) \right| \right) ^{%
\frac{1}{2}}; \\ 
\\ 
\max\limits_{1\leq i\leq n}\left| \alpha _{i}\right| \left(
\sum\limits_{i=1}^{n}\left| \alpha _{i}\right| ^{2\alpha }\right) ^{\frac{1}{%
2\alpha }}\left( \sum\limits_{i,j=1}^{n}\left| \left( z_{i},z_{j}\right)
\right| \right) ^{\frac{1}{2}}\left( \sum\limits_{i=1}^{n}\left(
\sum\limits_{j=1}^{n}\left| \left( z_{i},z_{j}\right) \right| \right)
^{\beta }\right) ^{\frac{1}{2\beta }}, \\ 
\hfill \ \ \text{if}\ \alpha >1,\ \frac{1}{\alpha }+\frac{1}{\beta }=1; \\ 
\left( \sum\limits_{i=1}^{n}\left| \alpha _{i}\right| ^{2\alpha }\right) ^{%
\frac{1}{2\alpha }}\left( \sum\limits_{i=1}^{n}\left| \alpha _{i}\right|
^{2\gamma }\right) ^{\frac{1}{2\gamma }}\left( \sum\limits_{i=1}^{n}\left(
\sum\limits_{j=1}^{n}\left| \left( z_{i},z_{j}\right) \right| \right)
^{\beta }\right) ^{\frac{1}{2\beta }} \\ 
\ \ \ \ \ \times \left( \sum\limits_{i=1}^{n}\left(
\sum\limits_{j=1}^{n}\left| \left( z_{i},z_{j}\right) \right| \right)
^{\delta }\right) ^{\frac{1}{2\delta }}\hfill \ \text{if}\ \alpha >1,\ \frac{%
1}{\alpha }+\frac{1}{\beta }=1\text{ and }\ \gamma >1,\ \frac{1}{\gamma }+%
\frac{1}{\delta }=1; \\ 
\\ 
\left( \sum\limits_{i=1}^{n}\left| \alpha _{i}\right| ^{2}\right) ^{\frac{1}{%
2}}\left( \sum\limits_{i=1}^{n}\left| \alpha _{i}\right| ^{2\alpha }\right)
^{\frac{1}{2\alpha }}\max\limits_{1\leq i\leq n}\left(
\sum\limits_{j=1}^{n}\left| \left( z_{i},z_{j}\right) \right| \right) ^{%
\frac{1}{2}}\left( \sum\limits_{i=1}^{n}\left( \sum\limits_{j=1}^{n}\left|
\left( z_{i},z_{j}\right) \right| \right) ^{\beta }\right) ^{\frac{1}{2\beta 
}}, \\ 
\hfill \ \ \text{if}\ \alpha >1,\ \frac{1}{\alpha }+\frac{1}{\beta }=1; \\ 
\\ 
\max\limits_{1\leq i\leq n}\left| \alpha _{i}\right| \left(
\sum\limits_{i=1}^{n}\left| \alpha _{i}\right| ^{2}\right) ^{\frac{1}{2}%
}\max\limits_{1\leq i\leq n}\left( \sum\limits_{j=1}^{n}\left| \left(
z_{i},z_{j}\right) \right| \right) ^{\frac{1}{2}}\left(
\sum\limits_{i,j=1}^{n}\left| \left( z_{i},z_{j}\right) \right| \right) ^{%
\frac{1}{2}}; \\ 
\\ 
\left( \sum\limits_{i=1}^{n}\left| \alpha _{i}\right| ^{2}\right) ^{\frac{1}{%
2}}\left( \sum\limits_{i=1}^{n}\left| \alpha _{i}\right| ^{2\gamma }\right)
^{\frac{1}{2\gamma }}\max\limits_{1\leq i\leq n}\left(
\sum\limits_{j=1}^{n}\left| \left( z_{i},z_{j}\right) \right| \right) ^{%
\frac{1}{2}}\left( \sum\limits_{i=1}^{n}\left( \sum\limits_{j=1}^{n}\left|
\left( z_{i},z_{j}\right) \right| \right) ^{\delta }\right) ^{\frac{1}{%
2\delta }}, \\ 
\hfill \ \ \text{if}\ \gamma >1,\ \frac{1}{\gamma }+\frac{1}{\delta }=1; \\ 
\sum\limits_{i=1}^{n}\left| \alpha _{i}\right| ^{2}\max\limits_{1\leq i\leq
n}\left( \sum\limits_{j=1}^{n}\left| \left( z_{i},z_{j}\right) \right|
\right) .
\end{array}
\right.
\end{equation*}
\end{corollary}

\section{Some Pe\v{c}ari\'{c} Type Inequalities}

We are now able to point out the following result which complements and
generalises the inequality (\ref{1.1} ) due to J. Pe\v{c}ari\'{c}.

\begin{theorem}
\label{t3.1}Let $x,y_{1},\dots ,y_{n}$ be vectors of an inner product space $%
\left( H;\left( \cdot ,\cdot \right) \right) $ and $c_{1},\dots ,c_{n}\in 
\mathbb{K}$. Then one has the inequalities: 
\begin{equation}
\left| \sum\limits_{i=1}^{n}c_{i}\left( x,y_{i}\right) \right| ^{2}
\label{3.1}
\end{equation}
\begin{equation*}
\leq \left\| x\right\| ^{2}\left( \sum\limits_{i=1}^{n}\left| c_{i}\right|
^{p}\left( \sum\limits_{j=1}^{n}\left| \left( y_{i},y_{j}\right) \right|
\right) \right) ^{\frac{1}{p}}\left( \sum\limits_{i=1}^{n}\left|
c_{i}\right| ^{q}\left( \sum\limits_{j=1}^{n}\left| \left(
y_{i},y_{j}\right) \right| \right) \right) ^{\frac{1}{q}}
\end{equation*}
\begin{equation*}
\leq \left\| x\right\| ^{2}\times \left\{ 
\begin{array}{l}
\max\limits_{1\leq i\leq n}\left| c_{i}\right|
^{2}\sum\limits_{i,j=1}^{n}\left| \left( y_{i},y_{j}\right) \right| ; \\ 
\\ 
\max\limits_{1\leq i\leq n}\left| c_{i}\right| \left(
\sum\limits_{i=1}^{n}\left| c_{i}\right| ^{\gamma q}\right) ^{\frac{1}{%
\gamma q}}\left( \sum\limits_{i,j=1}^{n}\left| \left( y_{i},y_{j}\right)
\right| \right) ^{\frac{1}{p}}\left( \sum\limits_{i=1}^{n}\left(
\sum\limits_{j=1}^{n}\left| \left( y_{i},y_{j}\right) \right| \right)
^{\delta }\right) ^{\frac{1}{\delta q}}, \\ 
\hfill \ \ \text{if}\ \gamma >1,\ \frac{1}{\gamma }+\frac{1}{\delta }=1; \\ 
\max\limits_{1\leq i\leq n}\left| c_{i}\right| \left(
\sum\limits_{i=1}^{n}\left| c_{i}\right| ^{q}\right) ^{\frac{1}{q}}\left(
\sum\limits_{i,j=1}^{n}\left| \left( y_{i},y_{j}\right) \right| \right) ^{%
\frac{1}{p}}\max\limits_{1\leq i\leq n}\left( \sum\limits_{j=1}^{n}\left|
\left( y_{i},y_{j}\right) \right| \right) ^{\frac{1}{q}}; \\ 
\\ 
\max\limits_{1\leq i\leq n}\left| c_{i}\right| \left(
\sum\limits_{i=1}^{n}\left| c_{i}\right| ^{\alpha p}\right) ^{\frac{1}{%
\alpha p}}\left( \sum\limits_{i,j=1}^{n}\left| \left( y_{i},y_{j}\right)
\right| \right) ^{\frac{1}{q}}\left( \sum\limits_{i=1}^{n}\left(
\sum\limits_{j=1}^{n}\left| \left( y_{i},y_{j}\right) \right| \right)
^{\beta }\right) ^{\frac{1}{p\beta }}, \\ 
\hfill \ \ \text{if}\ \alpha >1,\ \frac{1}{\alpha }+\frac{1}{\beta }=1; \\ 
\left( \sum\limits_{i=1}^{n}\left| c_{i}\right| ^{\alpha p}\right) ^{\frac{1%
}{\alpha p}}\left( \sum\limits_{i=1}^{n}\left| c_{i}\right| ^{\gamma
q}\right) ^{\frac{1}{\gamma q}}\left( \sum\limits_{i=1}^{n}\left(
\sum\limits_{j=1}^{n}\left| \left( y_{i},y_{j}\right) \right| \right)
^{\beta }\right) ^{\frac{1}{p\beta }} \\ 
\times \left( \sum\limits_{i=1}^{n}\left( \sum\limits_{j=1}^{n}\left| \left(
y_{i},y_{j}\right) \right| \right) ^{\delta }\right) ^{\frac{1}{\delta q}}\ 
\text{if}\ \alpha >1,\ \frac{1}{\alpha }+\frac{1}{\beta }=1\text{ and }%
\gamma >1,\ \frac{1}{\gamma }+\frac{1}{\delta }=1; \\ 
\\ 
\left( \sum\limits_{i=1}^{n}\left| c_{i}\right| ^{q}\right) ^{\frac{1}{q}%
}\left( \sum\limits_{i=1}^{n}\left| c_{i}\right| ^{\alpha p}\right) ^{\frac{1%
}{\alpha p}}\max\limits_{1\leq i\leq n}\left( \sum\limits_{j=1}^{n}\left|
\left( y_{i},y_{j}\right) \right| \right) ^{\frac{1}{q}} \\ 
\ \ \ \ \ \times \left( \sum\limits_{i=1}^{n}\left(
\sum\limits_{j=1}^{n}\left| \left( y_{i},y_{j}\right) \right| \right)
^{\beta }\right) ^{\frac{1}{p\beta }},\hfill \ \ \text{if}\ \alpha >1,\ 
\frac{1}{\alpha }+\frac{1}{\beta }=1; \\ 
\\ 
\max\limits_{1\leq i\leq n}\left| c_{i}\right| \left(
\sum\limits_{i=1}^{n}\left| c_{i}\right| ^{p}\right) ^{\frac{1}{p}%
}\max\limits_{1\leq i\leq n}\left( \sum\limits_{j=1}^{n}\left| \left(
y_{i},y_{j}\right) \right| \right) ^{\frac{1}{p}}\left(
\sum\limits_{i,j=1}^{n}\left| \left( y_{i},y_{j}\right) \right| \right) ^{%
\frac{1}{q}}; \\ 
\\ 
\left( \sum\limits_{i=1}^{n}\left| c_{i}\right| ^{p}\right) ^{\frac{1}{p}%
}\left( \sum\limits_{i=1}^{n}\left| c_{i}\right| ^{\gamma q}\right) ^{\frac{1%
}{\gamma q}}\max\limits_{1\leq i\leq n}\left( \sum\limits_{j=1}^{n}\left|
\left( y_{i},y_{j}\right) \right| \right) ^{\frac{1}{p}}\left(
\sum\limits_{i=1}^{n}\left( \sum\limits_{j=1}^{n}\left| \left(
y_{i},y_{j}\right) \right| \right) ^{\delta }\right) ^{\frac{1}{\delta q}},
\\ 
\ \ \ \ \ \ \ \hfill \ \ \text{if}\ \gamma >1,\ \frac{1}{\gamma }+\frac{1}{%
\delta }=1; \\ 
\left( \sum\limits_{i=1}^{n}\left| c_{i}\right| ^{p}\right) ^{\frac{1}{p}%
}\left( \sum\limits_{i=1}^{n}\left| c_{i}\right| ^{q}\right) ^{\frac{1}{q}%
}\max\limits_{1\leq i\leq n}\left( \sum\limits_{j=1}^{n}\left| \left(
y_{i},y_{j}\right) \right| \right) ;
\end{array}
\right.
\end{equation*}
where $p>1,$ $\frac{1}{p}+\frac{1}{q}=1.$
\end{theorem}

\begin{proof}
We note that 
\begin{equation*}
\sum\limits_{i=1}^{n}c_{i}\left( x,y_{i}\right) =\left(
x,\sum\limits_{i=1}^{n}\overline{c_{i}}y_{i}\right) .
\end{equation*}
Using Schwarz's inequality in inner product spaces, we have 
\begin{equation}
\left| \sum\limits_{i=1}^{n}c_{i}\left( x,y_{i}\right) \right| ^{2}\leq
\left\| x\right\| ^{2}\left\| \sum\limits_{i=1}^{n}\overline{c_{i}}%
y_{i}\right\| ^{2}.  \label{3.2}
\end{equation}
Finally, using Lemma \ref{l2.1} with $\alpha _{i}=\overline{c_{i}},$ $%
z_{i}=y_{i}$ $\left( i=1,\dots ,n\right) ,$ we deduce the desired inequality
(\ref{3.1}).
\end{proof}

\begin{remark}
\label{r3.2}If in (\ref{3.1}) we choose $p=q=2,$ we obtain amongst others,
the result $\left( \ref{1.1}\right) $ due to J. Pe\v{c}ari\'{c}.
\end{remark}

\section{Some Results of Bombieri Type}

The following results of Bombieri type hold.

\begin{theorem}
\label{t4.1}Let $x,y_{1},\dots ,y_{n}\in H.$ Then one has the inequality: 
\begin{multline}
\sum\limits_{i=1}^{n}\left| \left( x,y_{i}\right) \right| ^{2}  \label{4.1}
\\
\leq \left\| x\right\| \left[ \sum\limits_{i=1}^{n}\left| \left(
x,y_{i}\right) \right| ^{p}\left( \sum\limits_{j=1}^{n}\left| \left(
y_{i},y_{j}\right) \right| \right) \right] ^{\frac{1}{2p}} \\
\times \left[ \sum\limits_{i=1}^{n}\left| \left( x,y_{i}\right) \right|
^{q}\left( \sum\limits_{j=1}^{n}\left| \left( y_{i},y_{j}\right) \right|
\right) \right] ^{\frac{1}{2q}}
\end{multline}
\begin{equation*}
\leq \left\| x\right\| \times \left\{ 
\begin{array}{l}
\max\limits_{1\leq i\leq n}\left| \left( x,y_{i}\right) \right| \left(
\sum\limits_{i,j=1}^{n}\left| \left( y_{i},y_{j}\right) \right| \right) ^{%
\frac{1}{2}}; \\ 
\\ 
\max\limits_{1\leq i\leq n}\left| \left( x,y_{i}\right) \right| ^{\frac{1}{2}%
}\left( \sum\limits_{i=1}^{n}\left| \left( x,y_{i}\right) \right| ^{\gamma
q}\right) ^{\frac{1}{2\gamma q}}\left( \sum\limits_{i,j=1}^{n}\left| \left(
y_{i},y_{j}\right) \right| \right) ^{\frac{1}{2p}}\left(
\sum\limits_{i=1}^{n}\left( \sum\limits_{j=1}^{n}\left| \left(
y_{i},y_{j}\right) \right| \right) ^{\delta }\right) ^{\frac{1}{2\delta q}},
\\ 
\hfill \ \ \text{if}\ \gamma >1,\ \frac{1}{\gamma }+\frac{1}{\delta }=1; \\ 
\\ 
\max\limits_{1\leq i\leq n}\left| \left( x,y_{i}\right) \right| ^{\frac{1}{2}%
}\left( \sum\limits_{i=1}^{n}\left| \left( x,y_{i}\right) \right|
^{q}\right) ^{\frac{1}{2q}}\left( \sum\limits_{i,j=1}^{n}\left| \left(
y_{i},y_{j}\right) \right| \right) ^{\frac{1}{2p}}\max\limits_{1\leq i\leq
n}\left( \sum\limits_{j=1}^{n}\left| \left( y_{i},y_{j}\right) \right|
\right) ^{\frac{1}{2q}}; \\ 
\\ 
\max\limits_{1\leq i\leq n}\left| \left( x,y_{i}\right) \right| ^{\frac{1}{2}%
}\left( \sum\limits_{i=1}^{n}\left| \left( x,y_{i}\right) \right| ^{\alpha
p}\right) ^{\frac{1}{2\alpha \beta }}\left( \sum\limits_{i,j=1}^{n}\left|
\left( y_{i},y_{j}\right) \right| \right) ^{\frac{1}{2q}}\left(
\sum\limits_{i=1}^{n}\left( \sum\limits_{j=1}^{n}\left| \left(
y_{i},y_{j}\right) \right| \right) ^{\beta }\right) ^{\frac{1}{p\beta }}, \\ 
\hfill \ \ \text{if}\ \alpha >1,\ \frac{1}{\alpha }+\frac{1}{\beta }=1; \\ 
\left( \sum\limits_{i=1}^{n}\left| \left( x,y_{i}\right) \right| ^{\alpha
p}\right) ^{\frac{1}{2\alpha p}}\left( \sum\limits_{i=1}^{n}\left| \left(
x,y_{i}\right) \right| ^{\gamma q}\right) ^{\frac{1}{2\gamma q}}\left(
\sum\limits_{i=1}^{n}\left( \sum\limits_{j=1}^{n}\left| \left(
y_{i},y_{j}\right) \right| \right) ^{\beta }\right) ^{\frac{1}{2p\beta }} \\ 
\times \left( \sum\limits_{i=1}^{n}\left( \sum\limits_{j=1}^{n}\left| \left(
y_{i},y_{j}\right) \right| \right) ^{\delta }\right) ^{\frac{1}{2\delta q}%
}\hfill \ \text{if}\ \alpha >1,\ \frac{1}{\alpha }+\frac{1}{\beta }=1\text{
and }\ \gamma >1,\ \frac{1}{\gamma }+\frac{1}{\delta }=1; \\ 
\\ 
\left( \sum\limits_{i=1}^{n}\left| \left( x,y_{i}\right) \right| ^{q}\right)
^{\frac{1}{2q}}\left( \sum\limits_{i=1}^{n}\left| \left( x,y_{i}\right)
\right| ^{\alpha p}\right) ^{\frac{1}{2\alpha p}}\max\limits_{1\leq i\leq
n}\left( \sum\limits_{j=1}^{n}\left| \left( y_{i},y_{j}\right) \right|
\right) ^{\frac{1}{2p}} \\ 
\ \ \ \ \ \ \ \times \left( \sum\limits_{i=1}^{n}\left(
\sum\limits_{j=1}^{n}\left| \left( y_{i},y_{j}\right) \right| \right)
^{\beta }\right) ^{\frac{1}{2p\beta }},\hfill \ \ \text{if}\ \alpha >1,\ 
\frac{1}{\alpha }+\frac{1}{\beta }=1; \\ 
\\ 
\max\limits_{1\leq i\leq n}\left| \left( x,y_{i}\right) \right| ^{\frac{1}{2}%
}\left( \sum\limits_{i=1}^{n}\left| \left( x,y_{i}\right) \right|
^{p}\right) ^{\frac{1}{2p}}\max\limits_{1\leq i\leq n}\left(
\sum\limits_{j=1}^{n}\left| \left( y_{i},y_{j}\right) \right| \right) ^{%
\frac{1}{2p}}\left( \sum\limits_{i,j=1}^{n}\left| \left( y_{i},y_{j}\right)
\right| \right) ^{\frac{1}{2q}}; \\ 
\\ 
\left( \sum\limits_{i=1}^{n}\left| \left( x,y_{i}\right) \right| ^{p}\right)
^{\frac{1}{2p}}\left( \sum\limits_{i=1}^{n}\left| \left( x,y_{i}\right)
\right| ^{\gamma q}\right) ^{\frac{1}{2\gamma q}}\max\limits_{1\leq i\leq
n}\left( \sum\limits_{j=1}^{n}\left| \left( y_{i},y_{j}\right) \right|
\right) ^{\frac{1}{2p}} \\ 
\ \ \ \ \ \ \ \times \left( \sum\limits_{i=1}^{n}\left(
\sum\limits_{j=1}^{n}\left| \left( y_{i},y_{j}\right) \right| \right)
^{\delta }\right) ^{\frac{1}{2\delta q}},\hfill \ \ \text{if}\ \gamma >1,\ 
\frac{1}{\gamma }+\frac{1}{\delta }=1; \\ 
\left( \sum\limits_{i=1}^{n}\left| \left( x,y_{i}\right) \right| ^{p}\right)
^{\frac{1}{2p}}\left( \sum\limits_{i=1}^{n}\left| \left( x,y_{i}\right)
\right| ^{q}\right) ^{\frac{1}{2q}}\max\limits_{1\leq i\leq n}\left(
\sum\limits_{j=1}^{n}\left| \left( y_{i},y_{j}\right) \right| \right) ^{%
\frac{1}{2}},
\end{array}
\right.
\end{equation*}
where $p>1,\frac{1}{p}+\frac{1}{q}=1.$
\end{theorem}

\begin{proof}
The proof follows by Theorem \ref{t3.1} on choosing $c_{i}=\overline{\left(
x,y_{i}\right) },i\in \left\{ 1,...,n\right\} $ and taking the square root
in both sides of the inequalities involved. We omit the details.
\end{proof}

\begin{remark}
We observe, by the last inequality in $\left( \ref{4.1}\right) $, we get 
\begin{equation*}
\frac{\left( \sum\limits_{i=1}^{n}\left| \left( x,y_{i}\right) \right|
^{2}\right) ^{2}}{\left( \sum\limits_{i=1}^{n}\left| \left( x,y_{i}\right)
\right| ^{p}\right) ^{\frac{1}{p}}\left( \sum\limits_{i=1}^{n}\left| \left(
x,y_{i}\right) \right| ^{q}\right) ^{\frac{1}{q}}}\leq \left\| x\right\|
^{2}\max\limits_{1\leq i\leq n}\left( \sum\limits_{j=1}^{n}\left| \left(
y_{i},y_{j}\right) \right| \right) ,
\end{equation*}
where $p>1,\frac{1}{p}+\frac{1}{q}=1.$

If in this inequality we choose $p=q=2,$ then we recapture Bombieri's result 
$\left( \ref{1.3}\right) .$
\end{remark}

\end{document}